\documentclass[preprint,12pt]{elsarticle}

\usepackage[utf8]{inputenc}
\usepackage[T1]{fontenc}
\usepackage{newtxtext,newtxmath}
\usepackage{amsmath,bm}
\usepackage{graphicx}
\usepackage{booktabs,tabularx,array,multirow}
\usepackage{siunitx}
\usepackage{subcaption}
\usepackage{adjustbox}
\usepackage{xcolor}
\usepackage{placeins}
\usepackage{microtype}
\usepackage[hidelinks]{hyperref}
\usepackage[nameinlink,noabbrev]{cleveref}
\journal{Results in Engineering}
\biboptions{sort&compress}
\newcommand{\UF}{U_{\mathrm{F}}}

\begin{document}
\begin{frontmatter}

\title{Data-Driven Flutter Suppression via HODMD, DMDc and Constrained MPC}

\author[una]{Carlos Domingo Mendez Gaona\corref{cor1}}
\ead{cmendez@pol.una.py}
\cortext[cor1]{Corresponding author. ORCID: \href{https://orcid.org/0000-0003-3232-9188}{0000-0003-3232-9188}.}
\affiliation[una]{organization={GISPA Research Group, Facultad Politecnica, Universidad Nacional de Asuncion},
            city={San Lorenzo},
            country={Paraguay}}

\begin{abstract}
An interpretable data-to-control workflow is developed for
actuator-constrained flutter suppression by combining higher-order dynamic mode
decomposition (HODMD), dynamic mode decomposition with control (DMDc), and
model predictive control (MPC). The primary test case is a pitch--plunge
typical section with a trailing-edge flap and rational unsteady aerodynamics.
HODMD identifies the dominant frequency and growth rate in subcritical,
near-critical, and post-critical regimes from four measured structural
channels. Under 2\% RMS sensor noise, the near-critical growth-rate error is
$1.8\times10^{-3}$~s$^{-1}$, markedly lower than for standard DMD without
delay embedding. A five-state real reduced-order model is then formed from the
HODMD subspace, and its flap-input matrix is estimated from a small-amplitude
PRBS record using DMDc. On an independent chirp, the all-channel normalized
error is $2.4\times10^{-3}$. The resulting constrained MPC stabilizes the true
post-critical plant for an $8^\circ$ initial pitch perturbation with a
$1^\circ$ flap limit, whereas saturated LQR designed on the same model reaches
the prescribed validity cap. Three compact support studies assess
transferability: SU2 simulations show HODMD frequency recovery consistent with
the available spectral resolution; OpenFOAM dynamic-mesh calculations show
nearly linear articulated-flap moment authority over $2$--$8^\circ$ with
frequency-dependent gain; and a two-way SU2 Python-FSI case demonstrates
bounded closed-loop disturbance rejection with a PD baseline. These studies
support modal identification, actuator authority, and closed-loop feasibility,
but are not presented as CFD-level MPC validation. The workflow provides a
reproducible route from early-time measurements to constrained aeroelastic
control.
\end{abstract}

\end{frontmatter}

\section{Introduction}
Flutter remains a limiting instability for lightweight and increasingly
flexible aerospace structures. Active suppression has progressed from
analytical design and wind-tunnel demonstrations to flight programmes, but
model fidelity, actuator nonlinearities, robustness and certification remain
central obstacles~\cite{peloubet1984recent,mukhopadhyay2000benchmark,
livne2018aircraft,nasa2022flying}. Control-oriented aeroelastic models are
therefore required to preserve the unstable dynamics while remaining fast
enough for estimation and optimization.

Established reduced-order strategies include eigensystem realization from
response data~\cite{juang1985eigensystem}, Volterra representations for
nonlinear or CFD-generated aeroelastic responses~\cite{silva2005identification,
marzocca2004nonlinear,balajewicz2012reduced}, and impulse-response or
frequency-domain ROMs for benchmark configurations such as AGARD~445.6
\cite{silva2004development,yates1987agard}. These methods remain important
reference points. In parallel, dynamic mode decomposition (DMD) and Koopman
operator methods provide spectral descriptions directly from time-resolved
data~\cite{schmid2010dynamic,rowley2009spectral,tu2014dynamic,mezic2005spectral,
kutz2016dynamic}. DMD with control (DMDc) adds explicit actuation channels
\cite{proctor2016dynamic}, and Koopman/DMD predictors can be embedded in MPC
\cite{korda2018linear,kaiser2018sparse}.

Higher-order dynamic mode decomposition (HODMD) augments the measured data
with time delays, improving frequency and growth-rate identification when the
number of active dynamics exceeds the number of sensors or when measurements
are noisy~\cite{leclainche2017higher,vega2020higher,
leclainche2018analyzing}. The method has already shown value in flight-flutter
testing~\cite{leclainche2019newrobust,mendez2019aeroelastic,mendez2021new},
ground-vibration analysis~\cite{mendez2025gvt}, complex multiphase-flow
reconstruction~\cite{mendez2024bubble}, and wake-interaction analysis in
multibody configurations~\cite{beltran2019wake}. Recent developments include
low-memory on-the-fly HODMD~\cite{amor2023onthefly}, hierarchical feature
selection~\cite{corrochano2024hierarchical}, broader engineering
applications~\cite{groun2023beyond}, and data-driven flow-control sensitivity
analysis~\cite{lazpita2024sensibility}. More generally, machine-learning
methods are becoming an established component of aerospace modelling and
control~\cite{leclainche2023mlreview}.

The unresolved point addressed here is not whether HODMD can identify
flutter-related modes, nor whether MPC can suppress flutter in an analytical
model. Rather, it is whether an explicit and reproducible chain can be built
from a small set of measured channels to modal identification, input-enabled
reduced modelling, and actuator-constrained control, while clearly separating
the primary control result from supporting CFD evidence. Prior MPC studies
have shown the relevance of predictive constraints for two-dimensional
flap-driven wings~\cite{darabseh2022mpc}; the present work differs by deriving
the prediction model from response data through HODMD and DMDc and by testing
that model against a separate true plant.

The contributions are:
\begin{enumerate}
\item HODMD identification of the dominant stable and unstable aeroelastic
modes from four measurable structural channels, including a quantitative
comparison with standard DMD under 2\% RMS noise;
\item construction of a five-state real control-oriented ROM whose modal
structure comes from HODMD and whose flap-input matrix is estimated by DMDc,
with PRBS identification and independent chirp/free-response validation;
\item constrained MPC evaluation above the flutter speed under flap
saturation, including a large-initial-condition case in which saturated LQR
fails and a sensor-noise robustness test;
\item three deliberately separated support studies: SU2 modal validation,
OpenFOAM articulated-flap authority, and SU2 Python-FSI closed-loop
feasibility using bounded PD feedback.
\end{enumerate}

The scope is intentionally limited. The primary claim is early-time
suppression of the unstable mode of a linear benchmark. The OpenFOAM study
supports the physical plausibility of an articulated aerodynamic input, and
the SU2 Python-FSI study demonstrates that a bounded feedback load can be
closed around a coupled plant; neither is described as flap-based CFD flutter
suppression or as CFD-level MPC validation.

\section{Methodological workflow and scope}
\begin{figure}[t]
\centering
\includegraphics[width=\linewidth,height=.58\textheight,keepaspectratio]{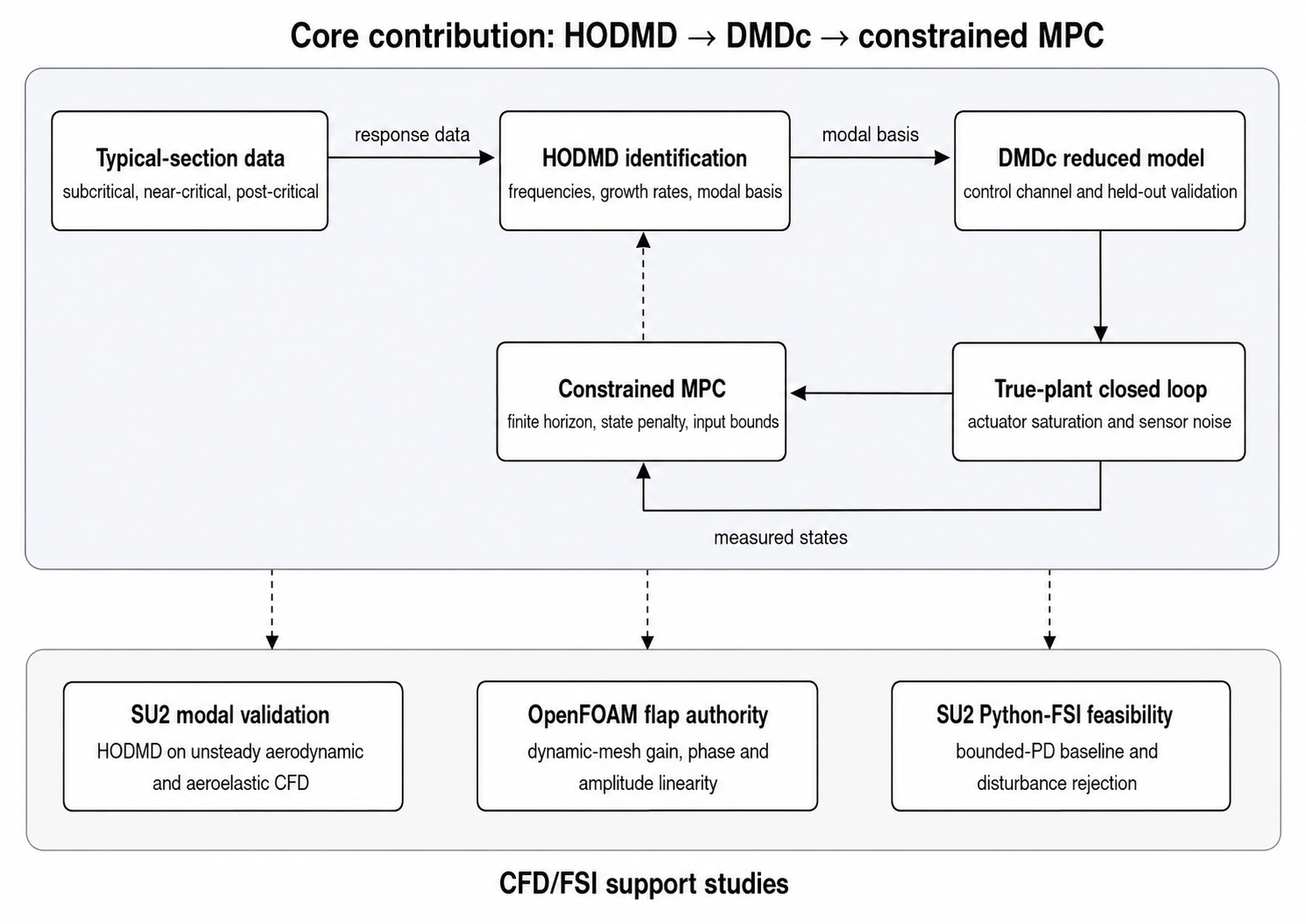}
\caption{Organization of the study. The primary result is the upper
HODMD--DMDc--MPC chain. SU2 and OpenFOAM are independent support studies and
do not constitute CFD-level validation of the MPC.}
\label{fig:workflow}
\end{figure}

Figure~\ref{fig:workflow} separates the paper's core contribution from the
three support studies. The upper chain is the only chain used to substantiate
the principal HODMD/DMDc/MPC claim. The lower studies answer narrower
engineering questions: whether HODMD remains useful on CFD data, whether a
moving articulated surface develops measurable unsteady authority, and
whether a feedback loop can be executed within a two-way Python-FSI driver.
This separation is maintained in the result statements and limitations.

\section{Aeroelastic benchmark, numerical verification and data generation}
\subsection{Structural model}
\begin{figure}[t]
\centering
\includegraphics[width=.92\linewidth,height=.50\textheight,keepaspectratio]{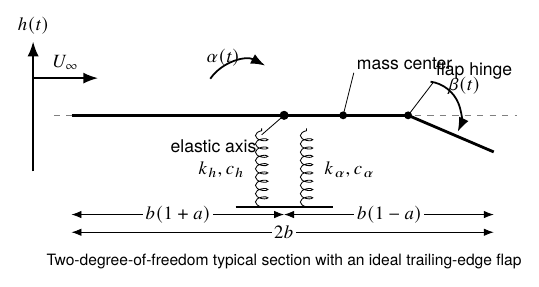}
\caption{Two-degree-of-freedom typical section used for the primary
HODMD/DMDc/MPC study. The structural coordinates are plunge $h$ and pitch
$\alpha$, and the control input is the ideal trailing-edge flap deflection
$\beta$.}
\label{fig:typical}
\end{figure}

Figure~\ref{fig:typical} defines the signs, characteristic dimensions and
control input used throughout the benchmark. The ideal flap is deliberately
kept separate from the later CFD assessment: it is the control input of the
reduced-order benchmark, whereas Section~\ref{sec:cfd-support} checks whether
an articulated surface exhibits a compatible aerodynamic authority.

The benchmark is a two-degree-of-freedom typical section, the classical
configuration of aeroelastic analysis~\cite{bisplinghoff1955aeroelasticity},
of semichord
$b$, with plunge $h$ (positive down) and pitch $\alpha$ (positive
nose-up) about an elastic axis located at $ab$ from midchord, and a
trailing-edge flap of deflection $\beta$ hinged at $cb$. Per unit span,
\begin{equation}
\begin{bmatrix} m & S_\alpha \\ S_\alpha & I_\alpha \end{bmatrix}
\begin{bmatrix} \ddot h \\ \ddot\alpha \end{bmatrix}
+
\begin{bmatrix} c_h & 0 \\ 0 & c_\alpha \end{bmatrix}
\begin{bmatrix} \dot h \\ \dot\alpha \end{bmatrix}
+
\begin{bmatrix} k_h & 0 \\ 0 & k_\alpha \end{bmatrix}
\begin{bmatrix} h \\ \alpha \end{bmatrix}
=
\begin{bmatrix} -L \\ M_{ea} \end{bmatrix},
\label{eq:eom}
\end{equation}
with static unbalance $S_\alpha = m x_\alpha b$ and uncoupled
frequencies $\omega_h=\sqrt{k_h/m}$, $\omega_\alpha=\sqrt{k_\alpha/I_\alpha}$.
The parameter values (Table~\ref{tab:params}) are illustrative choices
within textbook-typical nondimensional ranges; they do not represent any
specific aircraft.

\begin{table}[t]
\centering\small
\caption{Benchmark parameters (illustrative values; not representative
of any specific aircraft).}
\label{tab:params}
\begin{tabular}{lcc}
\toprule
Quantity & Symbol & Value \\
\midrule
Elastic-axis position & $a$ & $-0.4$ \\
Static unbalance & $x_\alpha$ & $0.2$ \\
Radius of gyration squared & $r_\alpha^2$ & $0.25$ \\
Mass ratio & $\mu$ & $20$ \\
Frequency ratio & $\omega_h/\omega_\alpha$ & $0.5$ \\
Structural damping ratios & $\zeta_h=\zeta_\alpha$ & $0.01$ \\
Flap hinge location & $c$ & $0.6$ \\
Semichord & $b$ & \SI{1.0}{m} \\
Air density & $\rho$ & \SI{1.225}{kg/m^3} \\
Pitch frequency & $\omega_\alpha$ & \SI{20}{rad/s} \\
\bottomrule
\end{tabular}
\end{table}

\subsection{Unsteady aerodynamics and state-space form}
\begin{figure}[t]
\centering
\includegraphics[width=\linewidth,height=.62\textheight,keepaspectratio]{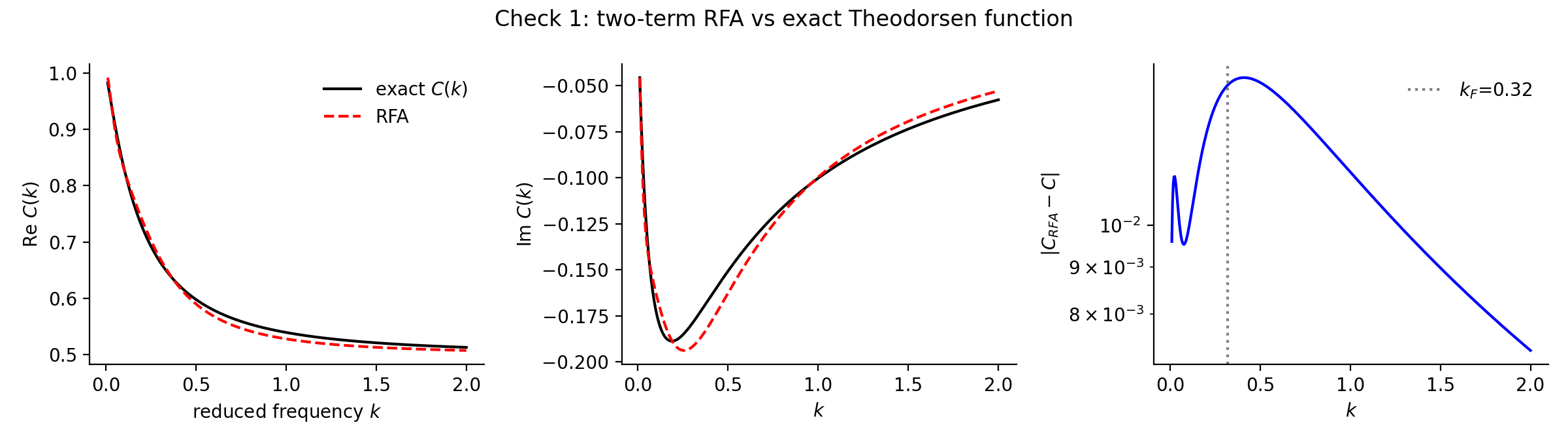}
\caption{Verification check~1: two-term rational approximation versus
the exact Theodorsen function. The dotted line marks the flutter reduced
frequency $k_{\mathrm F}=0.32$ of the benchmark.}
\label{fig:rfa}
\end{figure}

Lift and moment follow Theodorsen's framework~\cite{theodorsen1935general},
with noncirculatory
(apparent-mass) terms retained and the circulatory part driven by the
three-quarter-chord downwash
\begin{equation}
w_{3/4} = \dot h + U\alpha + b\left(\tfrac12-a\right)\dot\alpha
        + \frac{U T_{10}}{\pi}\,\beta ,
\label{eq:w34}
\end{equation}
where $T_{10}$ is the classical flap constant for hinge location $c$
(see the Appendix). In the core model the flap is ideal (massless,
commanded deflection, flap apparent-mass and $T_{11}\dot\beta$ terms
omitted); actuator dynamics are deferred to a sensitivity extension. For
time-domain simulation and control, the Theodorsen function $C(k)$ is
represented by a two-term rational (exponential) approximation of the
Wagner function of the type introduced by R.~T.
Jones~\cite{jones1940unsteady}, with the widely used coefficients $a_1=0.165$,
$b_1=0.0455$, $a_2=0.335$, $b_2=0.300$, which introduces two aerodynamic
lag states
\begin{equation}
\dot w_i = -\frac{b_i U}{b}\, w_i + w_{3/4}, \qquad i=1,2,
\label{eq:lags}
\end{equation}
and the circulatory lift
\begin{equation}
L_c = 2\pi\rho U b\left[(1-a_1-a_2)\,w_{3/4}
      + \frac{U}{b}\left(a_1 b_1 w_1 + a_2 b_2 w_2\right)\right].
\label{eq:Lc}
\end{equation}
The rational approximation is validated numerically against the exact
Theodorsen function (Hankel-function form): the maximum absolute error
is $0.0145$ for reduced frequencies $k\in[0.01,2]$, and $0.0143$ at the
flutter reduced frequency of this benchmark
(Fig.~\ref{fig:rfa}). With state vector
$\bm{x}=[h,\ \alpha,\ \dot h,\ \dot\alpha,\ w_1,\ w_2]^{\mathsf T}$ the
system is linear parameter-dependent in the flow speed $U$,
\begin{equation}
\dot{\bm{x}} = A(U)\,\bm{x} + B(U)\,\beta(t), \qquad
\bm{y} = C\bm{x} + \bm{\nu}(t),
\label{eq:ss}
\end{equation}
where $\bm{y}$ collects the four measurable channels
($h,\alpha,\dot h,\dot\alpha$) and $\bm{\nu}$ is optional measurement
noise. A quasi-steady baseline is obtained by setting $C(k)\equiv 1$
(no lag states).

\subsection{Flutter analysis}
\begin{figure}[t]
\centering
\includegraphics[width=\linewidth,height=.62\textheight,keepaspectratio]{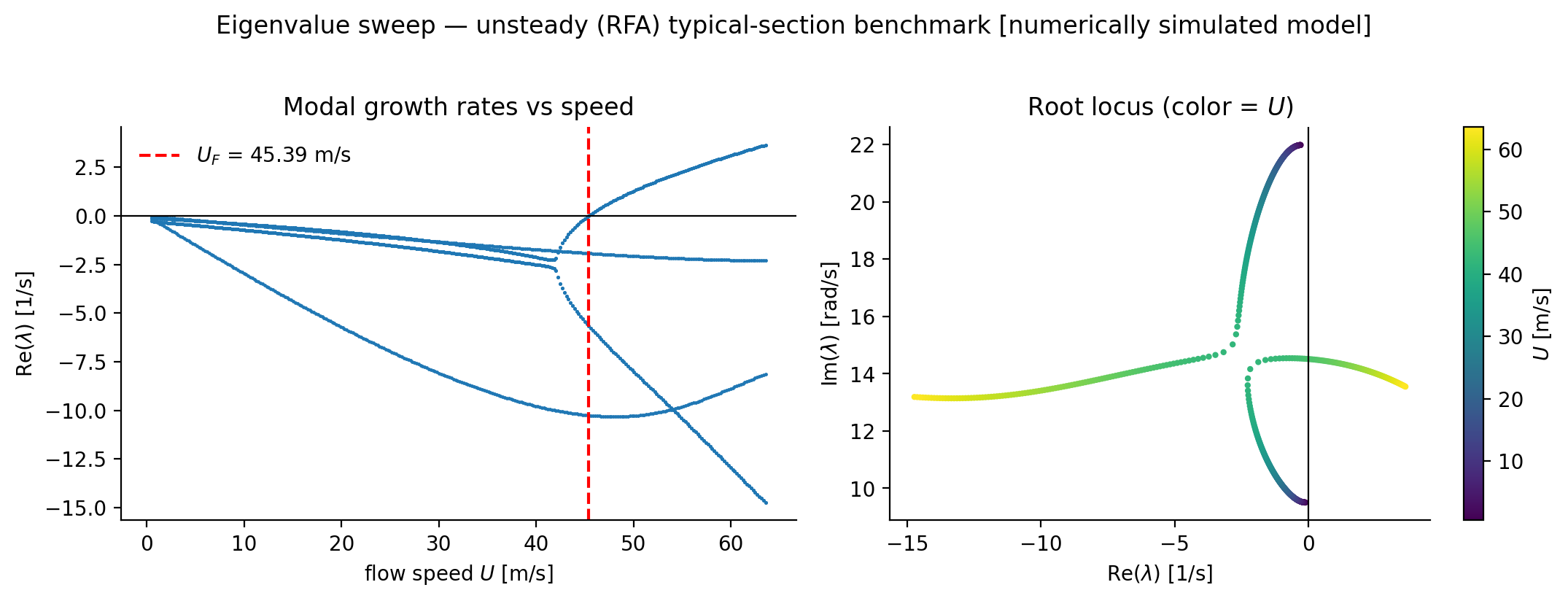}
\caption{Eigenvalue sweep of the unsteady benchmark: modal growth rates
versus flow speed (left) and root locus (right).}
\label{fig:sweep}
\end{figure}

An eigenvalue sweep of $A(U)$ (Fig.~\ref{fig:sweep}) locates the flutter
speed at $\UF=\SI{45.39}{m/s}$, with flutter frequency
\SI{2.31}{Hz} and reduced frequency $k_{\mathrm F}=0.32$; no static
divergence occurs within the swept range. The quasi-steady baseline
underpredicts the flutter speed by approximately \SI{34}{\percent}
($\UF^{\mathrm{qs}}=\SI{30.10}{m/s}$), quantifying the relevance of the
unsteady formulation for this configuration. Five verification checks
(rational-approximation accuracy; exact quasi-steady limit;
time-marching bracketing of $\UF$; integrator cross-check; structural
energy decay) pass and are documented with the released code.

\subsection{Data-generation and identification protocol}
All datasets are generated by exact zero-order-hold time marching of
Eq.~\eqref{eq:ss} (matrix exponential; cross-checked against a
tight-tolerance Runge--Kutta solution to $5.6\times10^{-10}$ relative
deviation) and are explicitly labeled as \emph{numerically simulated
data}. The time step resolves the fastest oscillatory mode with 50
samples per period; each record covers 45 periods of the dominant
flutter-related mode unless truncated. Post-critical records are
truncated when $|\alpha|$ exceeds \ang{15}, the stated validity cap of
the linear benchmark; all diverging simulations reported below are
truncated at this limit for the same reason.

\begin{figure}[t]
\centering
\includegraphics[width=\linewidth,height=.52\textheight,keepaspectratio]{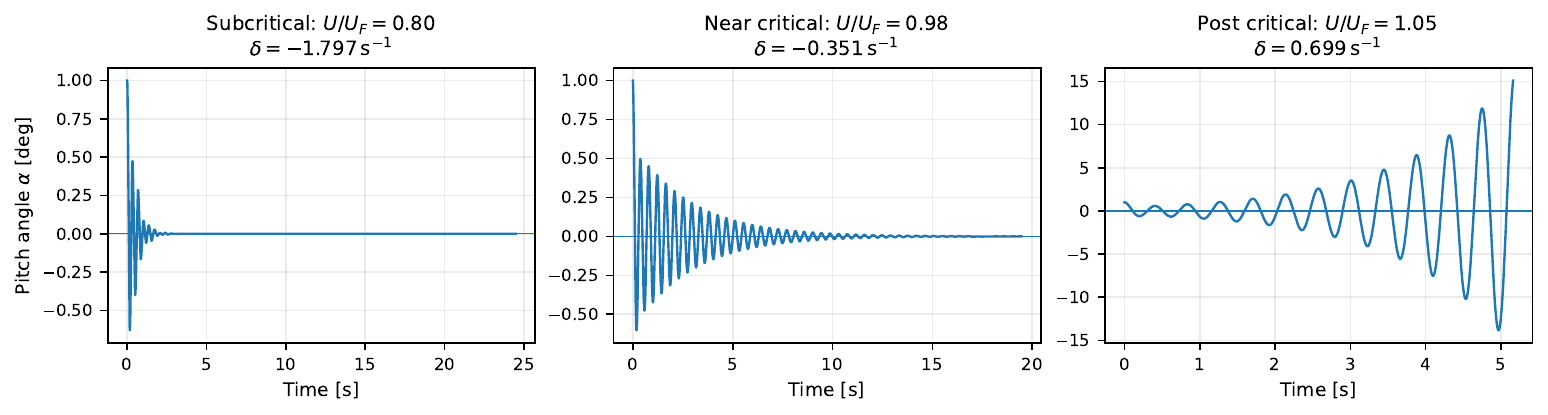}
\caption{Open-loop pitch response in the three identification regimes.
The post-critical record is truncated at $|\alpha|=\ang{15}$ and is used
only to identify the early-time unstable linear mode.}
\label{fig:openloop-regimes}
\end{figure}

\begin{figure}[t]
\centering
\includegraphics[width=\linewidth,height=.62\textheight,keepaspectratio]{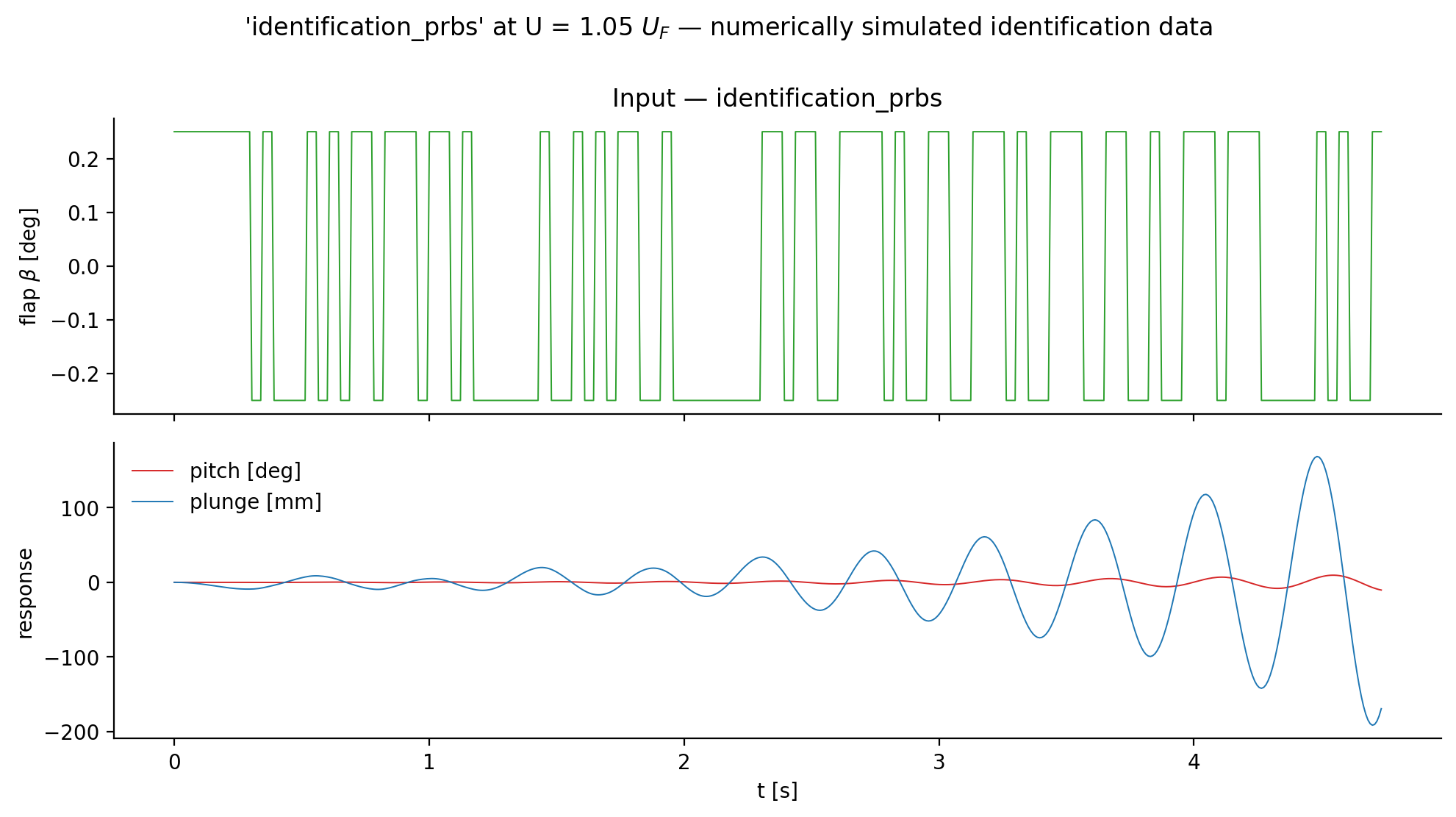}
\caption{Identification record: PRBS flap excitation
($\pm\ang{0.25}$) at $1.05\,\UF$ from zero initial conditions,
truncated at $|\alpha|=\ang{10}$. Numerically simulated identification
data.}
\label{fig:prbs}
\end{figure}

Three open-loop identification cases are generated from an initial
pitch perturbation $\alpha_0=\ang{1}$ at $0.80\,\UF$, $0.98\,\UF$, and
$1.05\,\UF$ (Table~\ref{tab:cases}; Fig.~\ref{fig:openloop-regimes}), plus a
robustness replica of the near-critical case with \SI{2}{\percent} RMS
additive Gaussian noise on the four measurable channels. For the
control-channel estimation, two additional records are generated at the
operating point $1.05\,\UF$ with zero initial conditions and small flap
excitation (\ang{0.25} amplitude, linear regime): a pseudo-random binary
sequence (PRBS, 43-ms hold) used for fitting, and a
\SIrange{0.5}{5}{Hz} chirp reserved exclusively for held-out validation
(Fig.~\ref{fig:prbs}).

\begin{table}[t]
\centering\small
\caption{Open-loop case matrix (numerically simulated data). $\sigma$
denotes the growth rate of the dominant oscillatory mode.}
\label{tab:cases}
\begin{tabular}{lcccccc}
\toprule
Case & $U/\UF$ & $f_s$ [Hz] & $N$ & $T$ [s] & $f_{\mathrm{dom}}$ [Hz]
 & $\sigma$ [1/s] \\
\midrule
Subcritical    & 0.80 & 145.5 & 3563 & 24.5 & 1.84 & $-1.797$ \\
Near-critical  & 0.98 & 115.6 & 2251 & 19.5 & 2.31 & $-0.351$ \\
Post-critical  & 1.05 & 114.9 & 594  & 5.2$^{\dagger}$ & 2.30 & $+0.699$ \\
Near-crit.\ + noise & 0.98 & 115.6 & 2251 & 19.5 & 2.31 & $-0.351$ \\
\bottomrule
\multicolumn{7}{l}{\footnotesize $^{\dagger}$Truncated at
$|\alpha|=\ang{15}$ (linear-model validity cap).}
\end{tabular}
\end{table}

\section{Early-time HODMD modal identification}
HODMD (DMD-$d$)~\cite{leclainche2017higher,vega2020higher} approximates
the measured snapshots by a finite expansion of spatial modes with
complex exponents,
\begin{equation}
\bm{y}(t_k) \approx \sum_{m=1}^{M} a_m\,\bm{\phi}_m\,
e^{(\delta_m + \mathrm{i}\omega_m)\,t_k},
\label{eq:hodmd}
\end{equation}
computed from $d$ time-lagged snapshot copies with two SVD-based
dimension reductions (tolerances $\varepsilon_1$) and amplitude-based
mode truncation. The sensor-realistic configuration uses only the four
measurable channels; the six-state data (including the aerodynamic lag
states) serve exclusively as a complementary internal check.

\begin{figure}[t]
\centering
\includegraphics[width=\linewidth,height=.62\textheight,keepaspectratio]{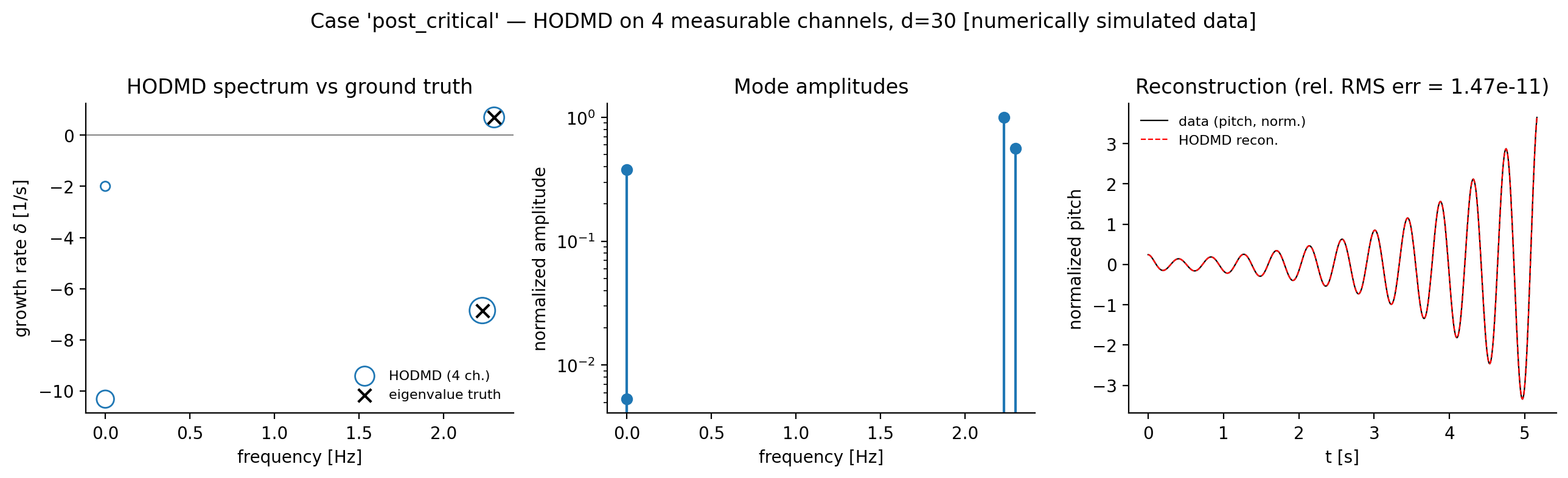}
\caption{HODMD on the short post-critical record (four measurable
channels): spectrum versus eigenvalue ground truth, mode amplitudes, and
reconstruction. Numerically simulated data.}
\label{fig:hodmdpost}
\end{figure}

On clean data, HODMD with $d=30$ recovers the dominant frequency and
growth rate of all three regimes to machine precision
(Table~\ref{tab:hodmd}), including the unstable mode
($\delta=+0.699$~1/s, \SI{2.298}{Hz}) from the short post-critical
record of approximately twelve periods (Fig.~\ref{fig:hodmdpost}). This
result is interpreted strictly as \emph{early-time identification of the
linear unstable mode}, not as a nonlinear post-flutter or limit-cycle
analysis. Secondary, more heavily damped modes are also recovered.

\begin{table}[t]
\centering\small
\caption{HODMD identification versus eigenvalue ground truth (dominant
mode, four measurable channels, $d=30$).}
\label{tab:hodmd}
\begin{tabular}{lccccc}
\toprule
Case & $f_{\mathrm{true}}$ [Hz] & $f_{\mathrm{est}}$ [Hz]
 & $\delta_{\mathrm{true}}$ [1/s] & $\delta_{\mathrm{est}}$ [1/s]
 & Recon.\ error \\
\midrule
Subcritical   & 1.8382 & 1.8382 & $-1.79677$ & $-1.79677$ & $1.3\times10^{-11}$ \\
Near-critical & 2.3130 & 2.3130 & $-0.35129$ & $-0.35129$ & $1.2\times10^{-11}$ \\
Post-critical & 2.2978 & 2.2978 & $+0.69899$ & $+0.69899$ & $1.5\times10^{-11}$ \\
Near-crit.\ + \SI{2}{\percent} noise
              & 2.3130 & 2.3127 & $-0.35129$ & $-0.34949$ & $2.1\times10^{-2}$ \\
\bottomrule
\end{tabular}
\end{table}

\begin{figure}[t]
\centering
\includegraphics[width=\linewidth,height=.62\textheight,keepaspectratio]{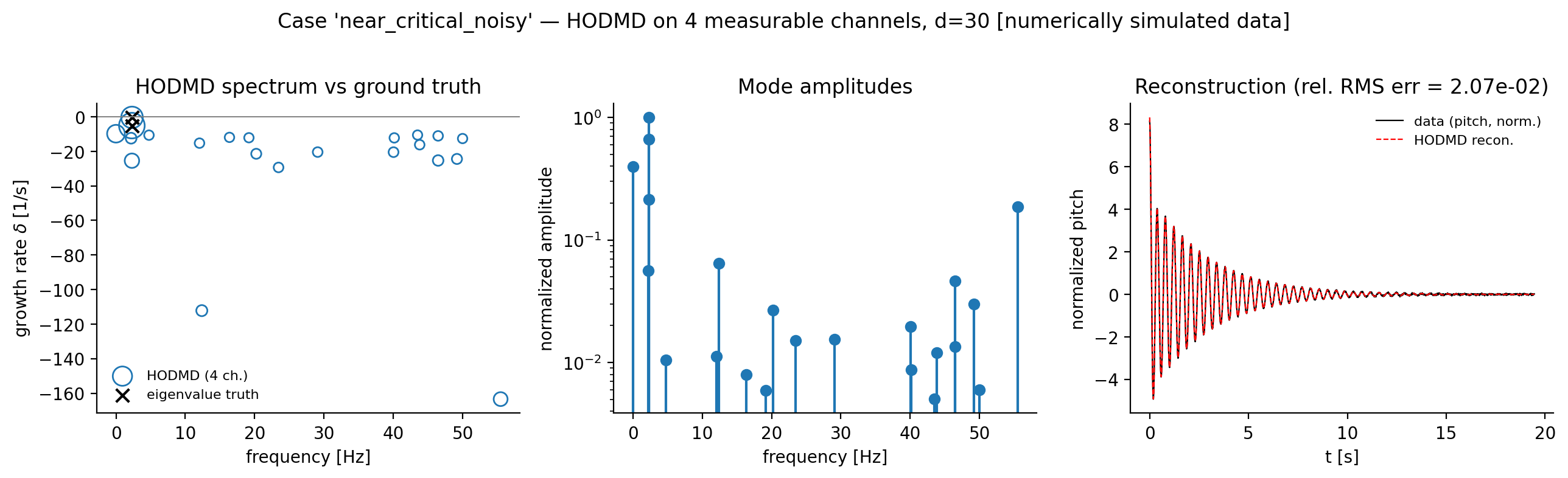}
\caption{HODMD on the near-critical record with \SI{2}{\percent} RMS
measurement noise: low-amplitude spurious modes appear alongside the
physical pairs; the reconstruction error saturates at the noise floor.}
\label{fig:hodmdnoisy}
\end{figure}

\begin{figure}[t]
\centering
\includegraphics[width=\linewidth,height=.62\textheight,keepaspectratio]{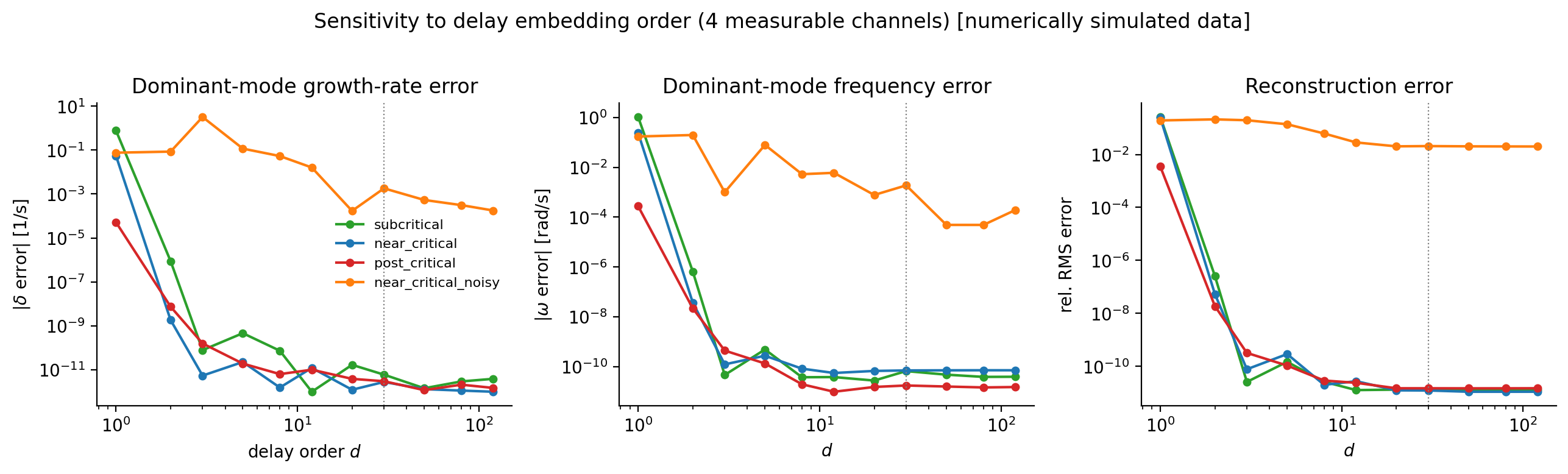}
\caption{Sensitivity of HODMD to the delay order $d$ (four measurable
channels): dominant-mode growth-rate and frequency errors and
reconstruction error. The dotted line marks $d=30$ used in the core
study.}
\label{fig:dsens}
\end{figure}

Under \SI{2}{\percent} RMS noise the dominant growth-rate error is
$1.8\times10^{-3}$~1/s (\SI{0.5}{\percent} relative) and the
reconstruction error saturates at the noise floor
($2.0\times10^{-2}$), indicating that the underlying clean dynamics are
reconstructed without overfitting the noise (Fig.~\ref{fig:hodmdnoisy}).
The delay order is decisive in this regime (Fig.~\ref{fig:dsens}):
standard DMD ($d=1$) yields a \SI{21}{\percent} relative growth-rate
error, decreasing monotonically below $10^{-3}$ for $d\gtrsim 20$. The
complementary six-state analysis improves the noisy-case error only
marginally ($1.2\times10^{-3}$ versus $1.8\times10^{-3}$), supporting
the interpretation that the delay embedding recovers the influence of
the unmeasured aerodynamic lag states without measuring them. A
practical caveat is reported honestly: with noise and intermediate $d$,
amplitude-based truncation retains numerous low-amplitude spurious
modes; physically meaningful modes are therefore isolated by combining
amplitude dominance with robustness of $(\delta,\omega)$ across delay
orders, the criterion adopted for the control-oriented model below.

\section{Control-oriented reduced model via HODMD and DMDc}
\subsection{Mode selection and real block realization}
Physically meaningful modes are retained if their amplitude exceeds
\SI{1}{\percent} of the maximum and their $(\delta,\omega)$ pair is
matched within $0.05$ across $d\in\{20,30,50\}$. Applied to the unforced
post-critical record, the criterion retains three modes
(Table~\ref{tab:modes}): the unstable pair, a damped pair, and one real
mode, giving a five-state model. Each complex-conjugate pair is realized
as a real $2\times2$ block,
\begin{equation}
\dot{\bm z}_m = \begin{bmatrix}\delta_m & -\omega_m\\
\omega_m & \delta_m\end{bmatrix}\bm z_m , \qquad
\bm{y}_m = \left[\,2\,\mathrm{Re}\,\bm\phi_m,\ -2\,\mathrm{Im}\,
\bm\phi_m\,\right]\bm z_m ,
\label{eq:blocks}
\end{equation}
and real modes as $1\times1$ blocks, yielding $(A_r, C_r)$ entirely from
the HODMD output.

\begin{table}[t]
\centering\small
\caption{Modes retained by the amplitude--robustness criterion
(unforced post-critical record, four measurable channels).}
\label{tab:modes}
\begin{tabular}{lccc}
\toprule
Mode & $\delta$ [1/s] & $f$ [Hz] & Relative amplitude \\
\midrule
Damped pair   & $-6.856$  & 2.228 & 1.000 \\
Unstable pair & $+0.699$  & 2.298 & 0.559 \\
Real mode     & $-10.321$ & 0     & 0.379 \\
\bottomrule
\end{tabular}
\end{table}

\subsection{Control-channel estimation and cross-check}
The input matrix $\hat B$ is estimated by least squares on the PRBS
record with $(A_r, C_r)$ fixed to the HODMD modal structure---a
DMDc-type regression~\cite{proctor2016dynamic} restricted to the
identified modal subspace (fit $0.9976$). As a consistency cross-check,
an unconstrained DMDc on the delay-embedded PRBS snapshots yields
eigenvalues within $1.25\times10^{-4}$ of the HODMD-retained ones.

\subsection{Held-out validation}
\begin{figure}[t]
\centering
\includegraphics[width=\linewidth,height=.62\textheight,keepaspectratio]{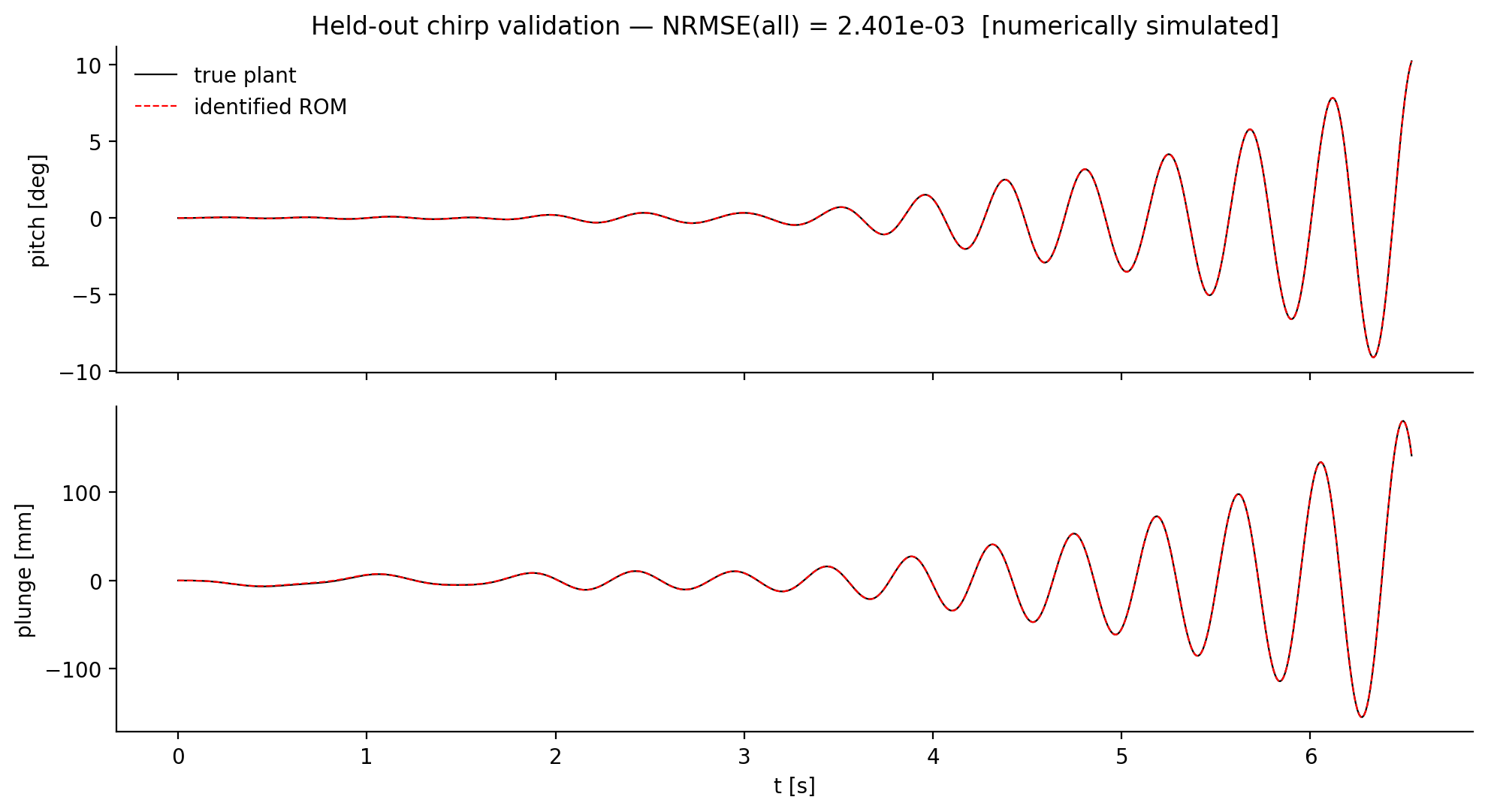}
\caption{Held-out validation: response of the identified five-state
model versus the true plant under the chirp input never used for
fitting.}
\label{fig:val}
\end{figure}

\begin{figure}[t]
\centering
\includegraphics[width=\linewidth,height=.62\textheight,keepaspectratio]{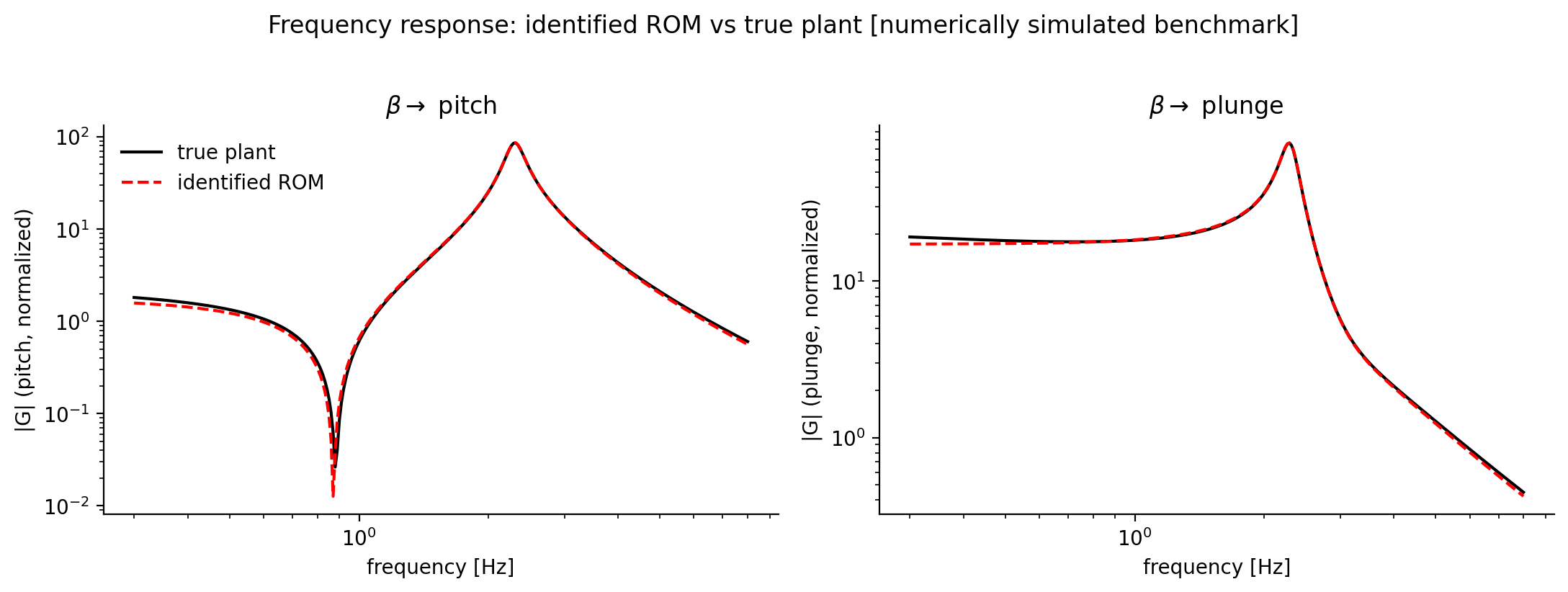}
\caption{Frequency response of the identified model versus the true
plant (flap to normalized outputs).}
\label{fig:freqresp}
\end{figure}

The identified five-state model is validated on the chirp record, which
was never used in any fitting step: it predicts the forced response of
the true six-state plant with normalized error $2.4\times10^{-3}$ over
all channels ($1.3\times10^{-3}$ on pitch; Fig.~\ref{fig:val}), and the
free response with error $7.7\times10^{-5}$. The frequency response of
the identified model matches that of the true plant across the
\SIrange{0.5}{5}{Hz} band (Fig.~\ref{fig:freqresp}). This comparison
doubles as a sanity check of the fully data-driven choice: the regressed
control channel is consistent with the analytically available one within
the benchmark. The fully data-driven pipeline is preferred for the main
study because it transfers to configurations where no analytical input
matrix exists; the hybrid alternative (HODMD dynamics with an analytical
input matrix) is retained only as supplementary verification.

\section{Constrained model predictive control}
The prediction model is the exact zero-order-hold discretization
$(A_d, B_d)$ of the identified $(A_r,\hat B)$ at the measurement rate
($\Delta t=\SI{8.7}{ms}$). The modal state is estimated by a
steady-state Kalman filter built from the identified $(A_d, C_r)$. At
each step the MPC solves
\begin{equation}
\min_{u_0,\dots,u_{N-1}}
\sum_{i=1}^{N}\bm y_i^{\mathsf T} Q\,\bm y_i + R\,u_{i-1}^2
\quad\text{s.t.}\quad
\bm z_{i+1}=A_d\bm z_i + B_d u_i,\quad |u_i|\le\beta_{\max},
\label{eq:mpc}
\end{equation}
with horizon $N=50$ (\SI{0.44}{s}), output weights
$Q=\mathrm{diag}(10,10,1,1)$ on the normalized measured channels, and
$R=50$. The condensed quadratic program with box constraints is solved
by a projected fast-gradient method. The baseline is an LQR designed on
the \emph{same} identified model and clipped to the \emph{same} flap
limits, so that the comparison isolates the effect of constraint
awareness. Both controllers act on the true six-state plant (i.e., in
the presence of model mismatch). The mean MPC solve time is
\SIrange{1.2}{1.3}{ms} per step against the \SI{8.7}{ms} sampling
interval; this figure is reported only as the computational cost in the
present software environment and does not constitute a claim of
real-time flight readiness.

\section{Closed-loop flutter-suppression results}
\subsection{Closed-loop suppression under actuator constraints}
\begin{figure}[t]
\centering
\includegraphics[width=\linewidth,height=.62\textheight,keepaspectratio]{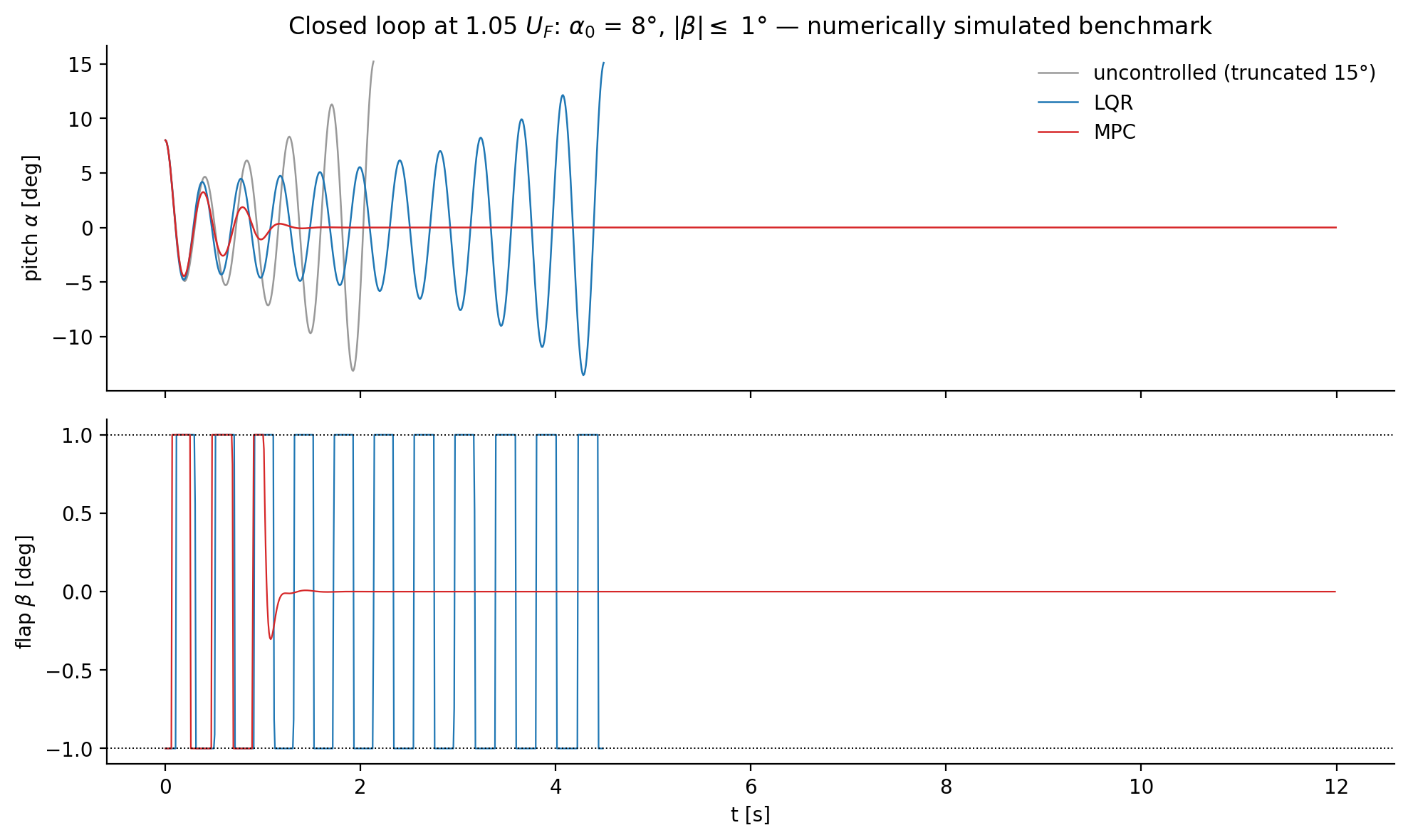}
\caption{Actuator-limited large-perturbation scenario
($\alpha_0=\ang{8}$, $|\beta|\le\ang{1}$, $1.05\,\UF$): the saturated
LQR diverges (truncated at \ang{15}) while the constraint-aware MPC
stabilizes the plant.}
\label{fig:cl}
\end{figure}

Table~\ref{tab:cl} summarizes the closed-loop scenarios at $1.05\,\UF$.
Three observations follow. First, when the flap constraint is inactive
or only weakly active (small perturbation, $\alpha_0=\ang{1}$), MPC and
LQR are indistinguishable, as expected since both minimize the same
quadratic criterion on the same identified model. Second, in the
actuator-limited large-perturbation case ($\alpha_0=\ang{8}$,
$\beta_{\max}=\ang{1}$; Fig.~\ref{fig:cl}), the saturated LQR fails to
stabilize the unstable plant---its response is reported truncated at
$|\alpha|=\ang{15}$, the validity limit of the linear
benchmark---whereas the constraint-aware MPC stabilizes the plant with
less than a quarter of the control energy expended by the diverging LQR.
This behavior is consistent with the known degradation of clipped linear
feedback on open-loop-unstable plants under sustained saturation and
constitutes, under the assumptions considered, the principal argument
for predictive control in actuator-limited flutter suppression. Third,
all stabilized responses correspond to early-time suppression of the
linear unstable mode; no statement is made about behavior beyond the
linear range.

\begin{table}[t]
\centering\footnotesize\setlength{\tabcolsep}{3.5pt}
\caption{Closed-loop metrics on the true plant at $1.05\,\UF$
(numerically simulated). Settling is to $|\alpha|<\ang{0.05}$. Diverging
runs are truncated at $|\alpha|=\ang{15}$ (linear-model validity cap).}
\label{tab:cl}
\begin{adjustbox}{max width=\linewidth}\begin{tabular}{llccccc}
\toprule
Scenario & Controller & Peak $|\alpha|$ [$^\circ$] & RMS $\alpha$
 [$^\circ$] & Settling [s] & $E_u$ [rad$^2\,$s] & Stable \\
\midrule
$\alpha_0{=}1^\circ$, $\beta{\le}5^\circ$ & LQR & 1.000 & 0.0632 & 0.279
 & $6.6\times10^{-5}$ & yes \\
$\alpha_0{=}1^\circ$, $\beta{\le}5^\circ$ & MPC & 1.000 & 0.0633 & 0.279
 & $6.5\times10^{-5}$ & yes \\
$\alpha_0{=}1^\circ$, $\beta{\le}1^\circ$ & LQR & 1.000 & 0.0689 & 0.418
 & $3.5\times10^{-5}$ & yes \\
$\alpha_0{=}1^\circ$, $\beta{\le}1^\circ$ & MPC & 1.000 & 0.0689 & 0.418
 & $3.5\times10^{-5}$ & yes \\
$\alpha_0{=}8^\circ$, $\beta{\le}1^\circ$ & LQR (sat.) &
 diverges$^{\dagger}$ & --- & --- & $1.4\times10^{-3}$ & no \\
$\alpha_0{=}8^\circ$, $\beta{\le}1^\circ$ & MPC & 8.000 & 0.791 & 1.471
 & $3.1\times10^{-4}$ & yes \\
Uncontrolled ($\alpha_0{=}1^\circ$ or $8^\circ$) & --- &
 diverges$^{\dagger}$ & --- & --- & 0 & no \\
\bottomrule
\multicolumn{7}{l}{\footnotesize $^{\dagger}$Truncated at
$|\alpha|=\ang{15}$ due to the validity limit of the linear benchmark.}
\end{tabular}
\end{adjustbox}
\end{table}

\subsection{Sensor-noise robustness}
\begin{figure}[t]
\centering
\includegraphics[width=\linewidth,height=.62\textheight,keepaspectratio]{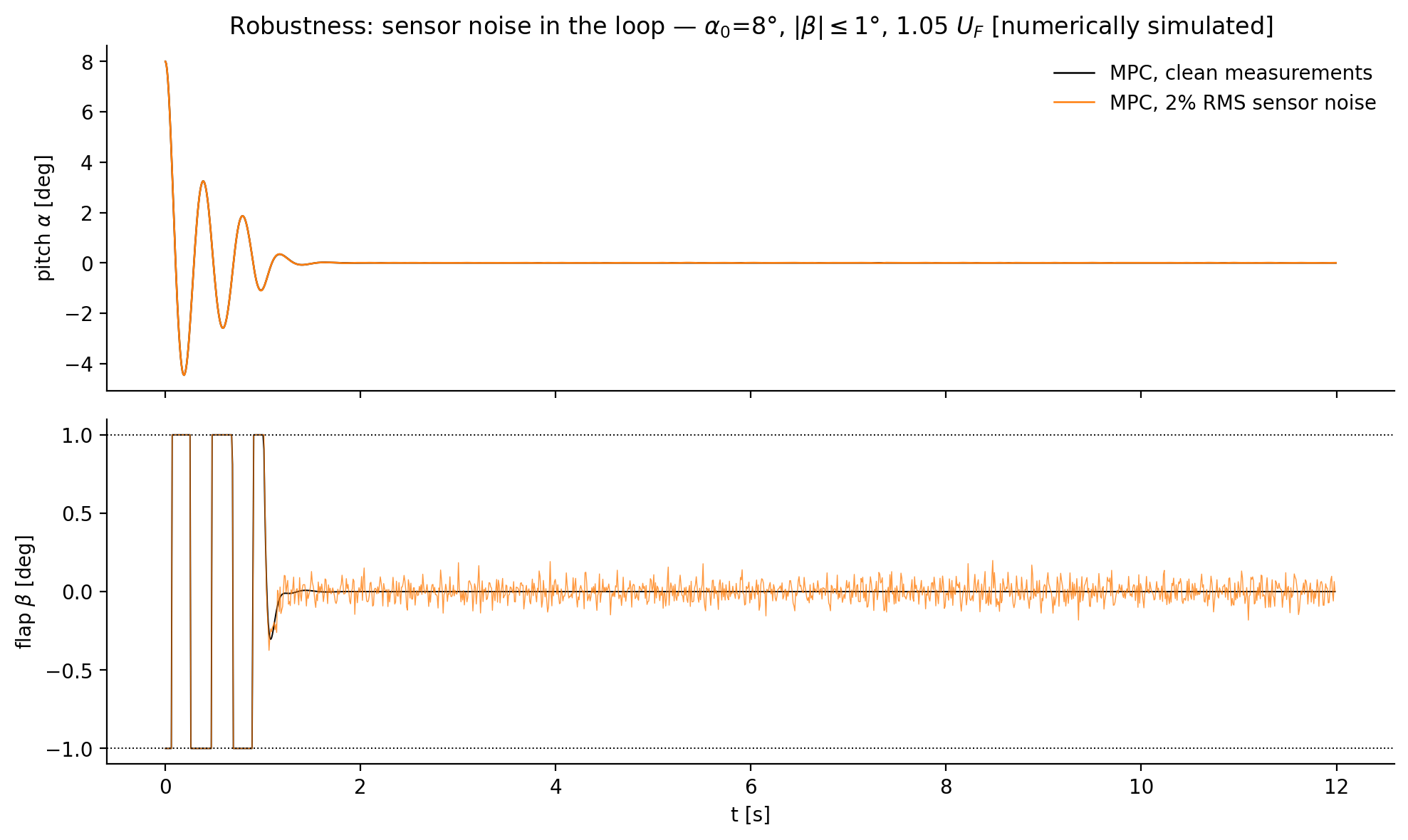}
\caption{Sensor-noise robustness ($\alpha_0=\ang{8}$,
$|\beta|\le\ang{1}$): MPC with clean versus \SI{2}{\percent}-RMS-noisy
measurements.}
\label{fig:noise}
\end{figure}

The most demanding stabilized scenario is repeated with Gaussian sensor
noise of \SI{2}{\percent} RMS per channel (referenced to the clean
closed-loop trajectory) injected into all four measurements, keeping the
same identified model and the same Kalman estimator
(Table~\ref{tab:noise}; Fig.~\ref{fig:noise}). The closed loop remains
stable and the transient metrics are unchanged; the noise manifests as a
small control jitter (\ang{0.059} RMS, well inside the \ang{1} limit)
and a residual pitch fluctuation an order of magnitude below the
settling threshold, at a \SI{3.9}{\percent} control-energy cost. The
saturated LQR remains unstable under noise. Within the limitations of
the proposed benchmark, this supports the claim that the pipeline
remains useful under sensor-realistic measurement conditions.

\begin{table}[t]
\centering\small
\caption{Clean versus noisy closed loop (MPC, $\alpha_0=\ang{8}$,
$\beta_{\max}=\ang{1}$). Steady-state statistics over $t>\SI{6}{s}$.}
\label{tab:noise}
\begin{tabular}{lcc}
\toprule
Metric & Clean & \SI{2}{\percent} RMS noise \\
\midrule
Peak $|\alpha|$ [$^\circ$] & 8.000 & 8.000 \\
RMS $\alpha$ [$^\circ$] & 0.7909 & 0.7909 \\
Settling to \ang{0.05} [s] & 1.471 & 1.471 \\
Control energy [rad$^2\,$s] & $3.09\times10^{-4}$
 & $3.21\times10^{-4}$ ($+3.9\%$) \\
Steady-state pitch jitter [$^\circ$ RMS] & $\sim 7\times10^{-11}$
 & $3.1\times10^{-3}$ (max $0.010$) \\
Steady-state control jitter [$^\circ$ RMS] & $\sim 6\times10^{-10}$
 & $5.9\times10^{-2}$ \\
Stable & yes & yes \\
\bottomrule
\end{tabular}
\end{table}

\subsection{What worked and what required care}
The elements that performed as intended were HODMD as an identification
engine across all regimes, including short-record unstable-mode
identification; the delay embedding compensating unmeasured aerodynamic
states; the modal-subspace control-channel regression, validated on
genuinely held-out data; constraint-aware MPC in the regime where
saturated linear feedback fails; and the noise robustness of the full
loop. The elements requiring care were the retention of spurious modes
by amplitude-only truncation under noise, resolved by the cross-$d$
robustness criterion; the \SI{34}{\percent} flutter-speed error of the
quasi-steady formulation, resolved by the rational approximation; and
the reporting of diverging runs, resolved by explicit truncation at the
linear-validity cap.


\section{CFD and FSI support studies}\label{sec:cfd-support}
The support studies were added after the reduced-order benchmark had been
verified. They are intentionally compact and answer specific questions rather
than forming a second control paper. SU2 is cited as the high-fidelity
multiphysics platform~\cite{economon2016su2}; the dynamic-mesh calculations
use the OpenFOAM finite-volume framework~\cite{weller1998openfoam}.
Table~\ref{tab:cfd-resolution} reports the minimum numerical descriptors needed
to judge these support calculations. Because these cases are not the primary
validation target, no Richardson extrapolation or formal grid-convergence
study was performed. The resulting claims are therefore restricted to
frequency consistency, harmonic trends, and closed-loop feasibility on the
reported meshes. The laminar and inviscid cases do not use wall functions, so
$y^+$ is not an applicable acceptance metric.

\begin{table}[t]
\centering
\caption{Numerical settings for the CFD/FSI support studies.}
\label{tab:cfd-resolution}
\footnotesize
\begin{tabularx}{\linewidth}{@{}lXXX@{}}
\toprule
Case & Spatial resolution & Flow model & Time and coupling settings \\
\midrule
SU2 plunging NACA0012 & 5233 points; 10216 volume elements &
Laminar Navier--Stokes; $M=0.30$; $Re=1000$ &
$\Delta t=2.3555\times10^{-3}$~s; 2000 steps; second-order dual time; CFL 1 \\
SU2 NACA64A010 & 6532 points; 9313 volume elements &
Euler; $M=0.85$; aeroelastic surface motion &
$\Delta t=1.7453\times10^{-3}$~s; 2000 steps; second-order dual time; CFL 100 \\
OpenFOAM articulated flap & 6828 hexahedra; 14098 points &
2-D laminar incompressible; $Re=100$ &
$\Delta t=2.0\times10^{-3}$~s; three cycles; $Co_{\max}=0.1015$ \\
SU2 Python-FSI & Fluid: 2370 points/2240 cells; structure: 150 points/116 elements &
Laminar Navier--Stokes; $M=0.05$; $Re=100$; two-way FSI &
$\Delta t=10^{-3}$~s; 4000 or 8000 steps; 30 outer iterations \\
\bottomrule
\end{tabularx}
\end{table}

The OpenFOAM mesh passed \texttt{checkMesh}, with maximum
non-orthogonality $13.65^\circ$, maximum skewness 0.677, and maximum aspect
ratio 2.60. The SU2 meshes were the distributed benchmark meshes used by the
corresponding test cases. None of the support calculations used wall
functions: the NACA64A010 case is inviscid and the remaining viscous cases are
laminar, so $y^+$ is not used as a turbulence-resolution criterion.

\subsection{SU2 modal validation}
Two unsteady SU2 cases were analysed using the same HODMD logic as the
benchmark. The prescribed plunging NACA0012 record contained 2000 samples at
$\Delta t=2.3555\times10^{-3}$~s, giving a record length of 4.71~s and an FFT
bin width of approximately $\Delta f=0.21$~Hz. The FFT lift peak was
16.87~Hz, while the continuous-time HODMD fit gave 16.98~Hz. Their 0.11~Hz
difference is smaller than one FFT bin and is therefore reported as
consistency within the available spectral resolution, rather than as a
sub-bin FFT comparison. The HODMD reconstruction had NRMSE 0.0072 and a
nearly neutral fitted growth rate of $-9.2\times10^{-5}$~s$^{-1}$. In the
coupled NACA64A010 aeroelastic case, HODMD identified 18.85~Hz versus an FFT
peak of 18.77~Hz, with a fitted growth rate of
$-5.0\times10^{-3}$~s$^{-1}$ and approximately 1.4\% reconstruction error
across lift, moment, plunge, and pitch. These cases support transfer of the
modal-identification step to CFD-generated data; they are not used to infer a
flutter boundary.

\begin{figure}[t]
\centering
\begin{subfigure}[t]{.49\linewidth}
\centering
\includegraphics[width=\linewidth,height=.34\textheight,keepaspectratio]{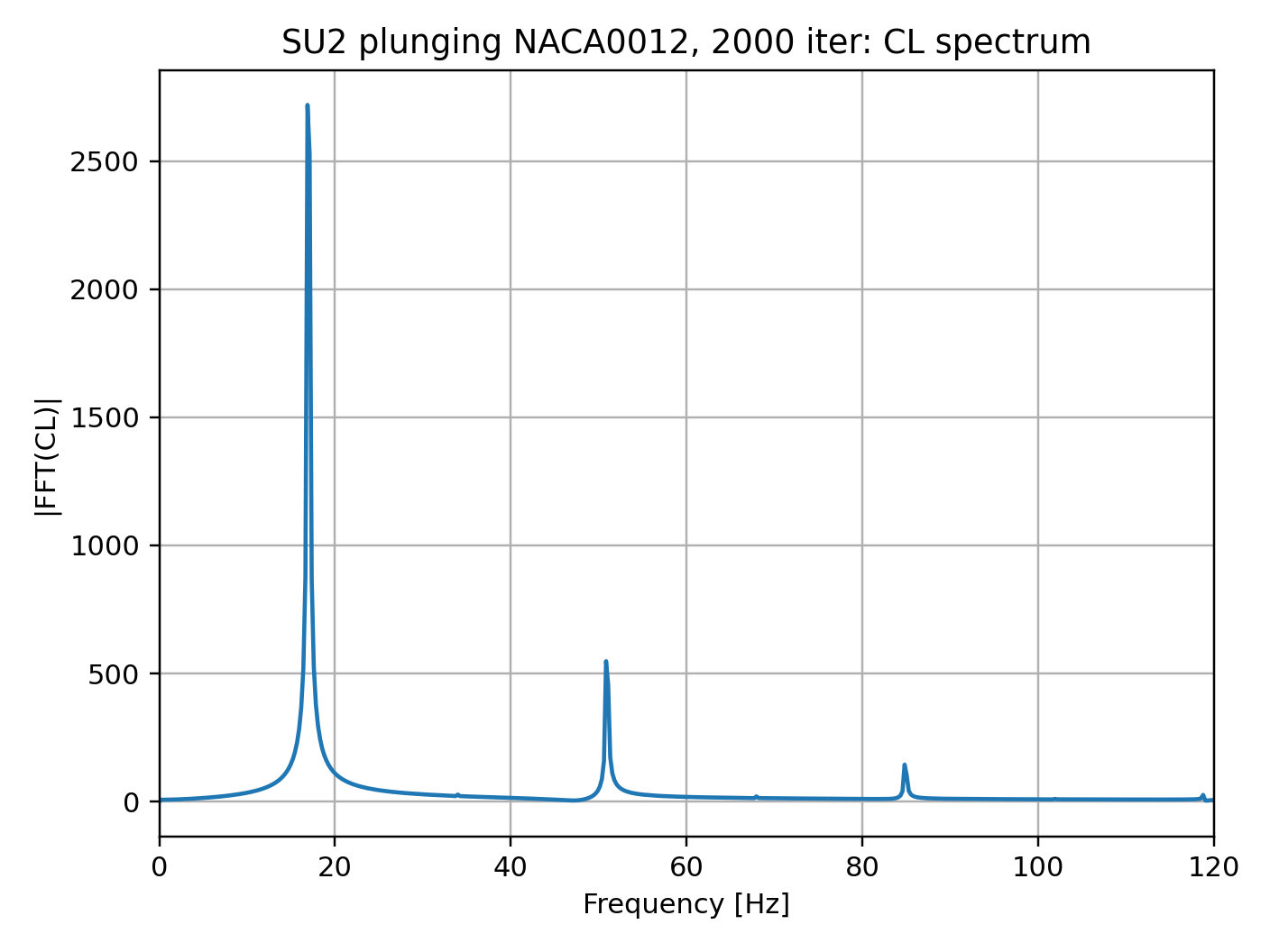}
\caption{Plunging NACA0012: FFT peak.}
\end{subfigure}\hfill
\begin{subfigure}[t]{.49\linewidth}
\centering
\includegraphics[width=\linewidth,height=.34\textheight,keepaspectratio]{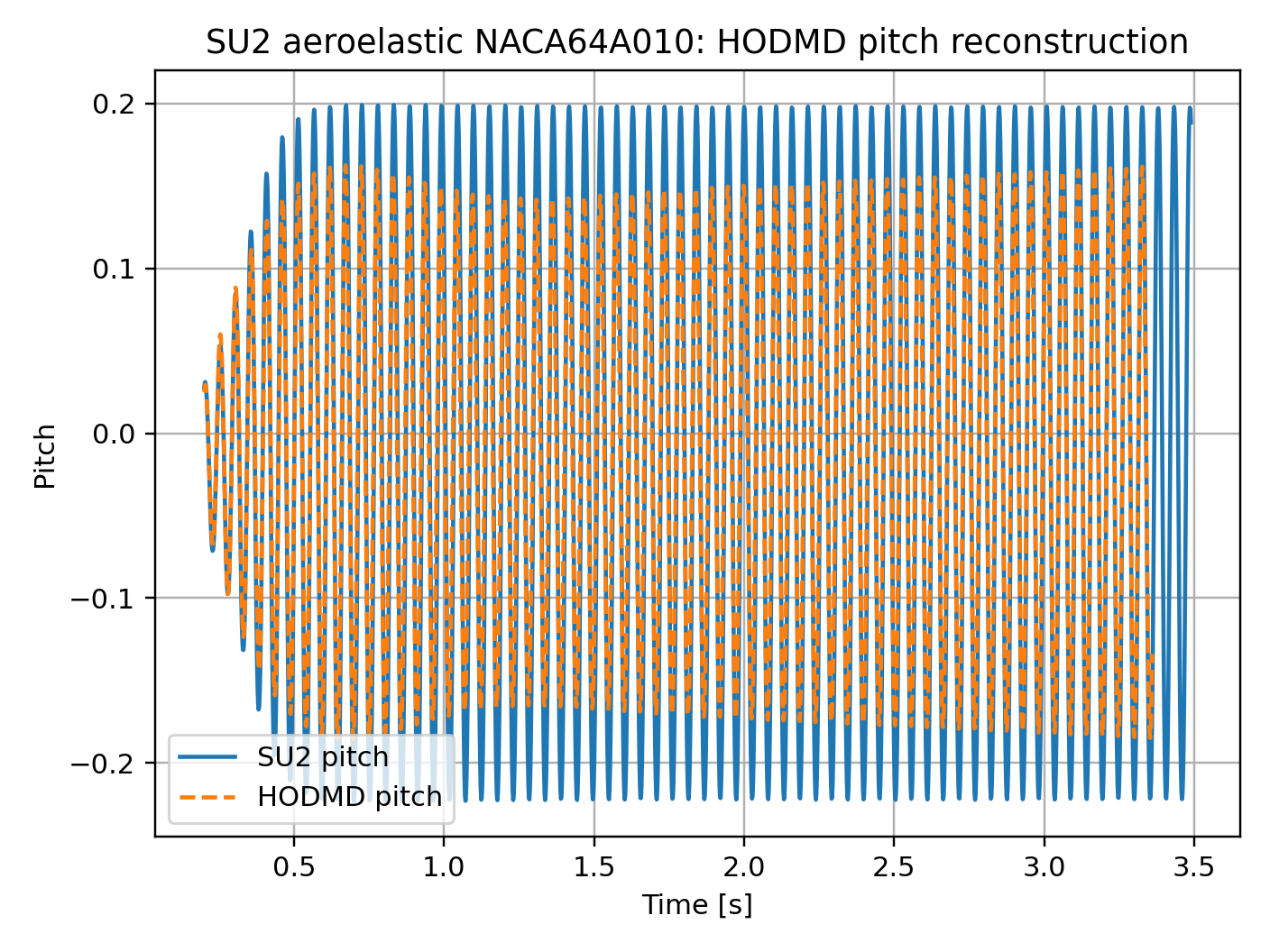}
\caption{NACA64A010: pitch reconstruction.}
\end{subfigure}
\caption{SU2 modal checks. HODMD recovers the dominant frequency in a
prescribed unsteady-aerodynamic case and reconstructs the coupled pitch
response in an aeroelastic case.}
\label{fig:su2-modal}
\end{figure}

\subsection{OpenFOAM articulated-flap authority}
A separate two-dimensional, laminar, low-Reynolds-number demonstrator was
constructed with a fixed main body, a rear articulated surface, and a
deforming mesh. The purpose was not to reproduce the typical-section
geometry, but to test the input assumption $\delta(t)\mapsto
(F_y,M_z)$ under an explicitly moving boundary. Forces and moments were
integrated separately on the main body and flap and then combined about a
common reference point.

\begin{figure}[t]
\centering
\includegraphics[width=.90\linewidth,height=.48\textheight,keepaspectratio]{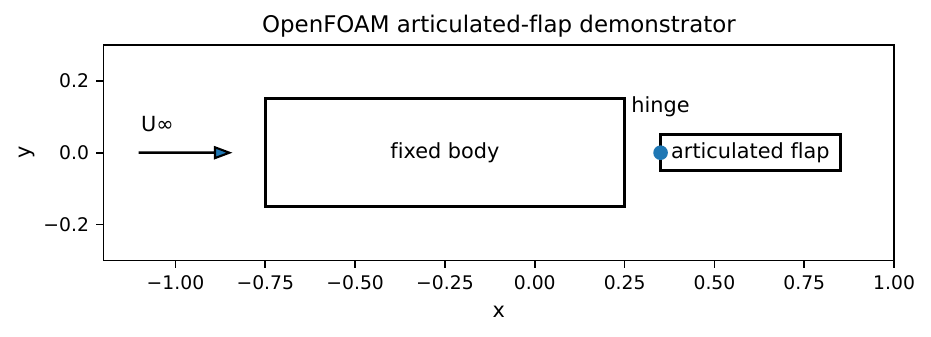}
\caption{OpenFOAM articulated-surface demonstrator used to assess unsteady
control authority. The simplified two-dimensional geometry, artificial hinge
gap and $Re=100$ flow are stated limitations; the result is used only as an
actuation-authority check.}
\label{fig:foam-geometry}
\end{figure}

The reported authority is dimensional. It is defined as
$G_{M\delta}=|M_{z,1}|/\delta_0$ in N\,m\,rad$^{-1}$ and
$G_{F\delta}=|F_{y,1}|/\delta_0$ in N\,rad$^{-1}$, with no dynamic-pressure
or chord normalization. The calculation uses $\rho=1$~kg\,m$^{-3}$,
$U_\infty=1$~m\,s$^{-1}$, a 0.02~m extruded depth, and moments about
$C_{\mathrm{ofR}}=(0,0,0.01)$~m. For $\delta_0=2^\circ$, $5^\circ$, and
$8^\circ$ at $\omega=\pi$~rad/s, $G_{M\delta}$ was 0.01011, 0.01009, and
0.01008~N\,m\,rad$^{-1}$, respectively: a variation of 0.23\%. The
lateral-force gain varied by 0.28\%, moment fits retained $R^2>0.995$, and
the second harmonic remained near 1.3\%. The amplitude study therefore
supports a provisional linear range of $\pm8^\circ$ for this demonstrator.
The frequency sweep, however, shows that a static gain is insufficient:
$G_{M\delta}$ rises from 0.00417 at $\pi/2$~rad/s to
0.03806~N\,m\,rad$^{-1}$ at $2\pi$~rad/s, while the phase changes from
$-60.0^\circ$ to $-22.5^\circ$.

\begin{figure}[t]
\centering
\begin{subfigure}[t]{.49\linewidth}
\centering
\includegraphics[width=\linewidth,height=.34\textheight,keepaspectratio]{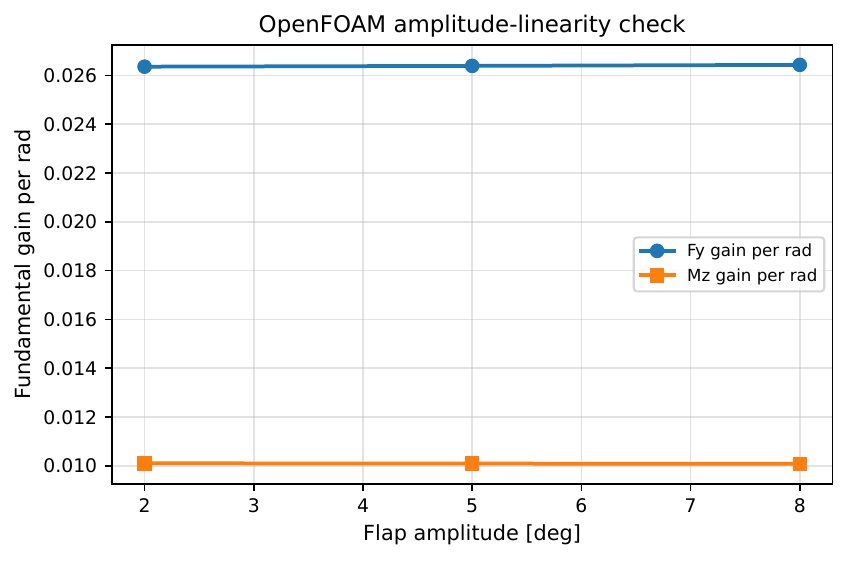}
\caption{Amplitude linearity.}
\end{subfigure}\hfill
\begin{subfigure}[t]{.49\linewidth}
\centering
\includegraphics[width=\linewidth,height=.34\textheight,keepaspectratio]{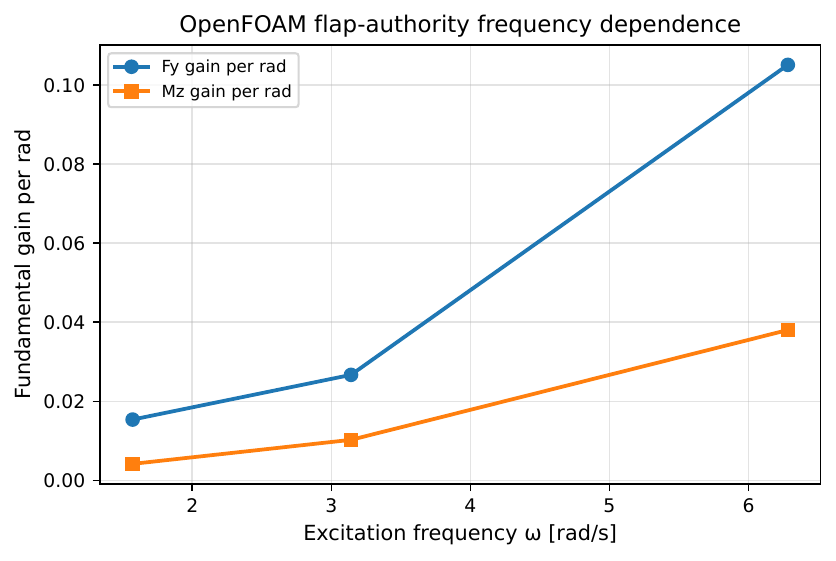}
\caption{Frequency dependence.}
\end{subfigure}
\caption{Dynamic-mesh articulated-flap authority. The nearly constant gain
with amplitude supports a locally linear actuation model, whereas the strong
frequency dependence motivates a dynamic input representation.}
\label{fig:foam-gains}
\end{figure}

\subsection{SU2 Python-FSI closed-loop feasibility}
\begin{figure}[t]
\centering
\includegraphics[width=\linewidth,height=.50\textheight,keepaspectratio]{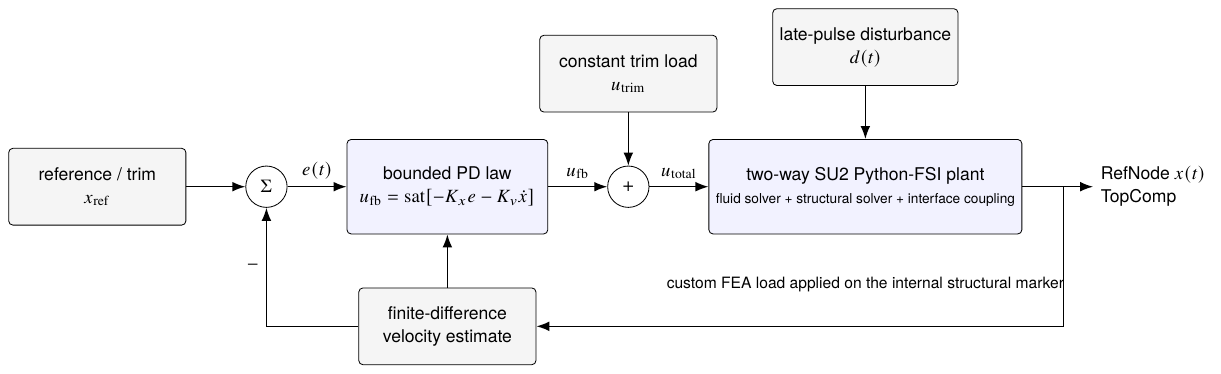}
\caption{Closed-loop architecture used only in the SU2 Python-FSI support
study. The feedback law is a bounded PD controller acting through a custom
structural load; it is not the MPC of the primary benchmark.}
\label{fig:fsi-loop}
\end{figure}

The final support study uses a two-way SU2 Python-FSI flexible-wall benchmark.
It does not contain the typical-section flap and therefore does not validate
the primary MPC. Its purpose is to verify the mechanics of a high-fidelity
closed loop: reading a structural response, computing a bounded control input,
applying a custom FEA load, and continuing the coupled fluid/structural solve.
Figure~\ref{fig:fsi-loop} shows the implemented baseline.

For the softened case, the dimensional Young's modulus was
$E=6.5\times10^{4}$~Pa and the structural density was 50~kg\,m$^{-3}$.
The coupled calculations used $\Delta t=10^{-3}$~s, 30 outer coupling
iterations, and a static relaxation factor of 0.1. The open-loop run
terminated at 3.303~s, whereas both a constant trim load and bounded feedback
completed the 4000-step horizon. A fair constant-load baseline showed that
static bias was sufficient to prevent divergence, so the feedback claim is
restricted to disturbance rejection. With $K_x=0.18$, $K_v=0.020$ and
$u_{\max}=0.0020$, a late pulse of $A_2=8\times10^{-4}$ produced a 63.1\%
reduction in peak RefNode deviation and a 74.6\% reduction in post-pulse RMS
relative to constant trim. Across $A_2=4\times10^{-4}$,
$8\times10^{-4}$, and $1.2\times10^{-3}$, the RMS reductions were 66.4\%,
74.6\%, and 63.7\%, respectively.

\begin{figure}[t]
\centering
\includegraphics[width=.90\linewidth,height=.46\textheight,keepaspectratio]{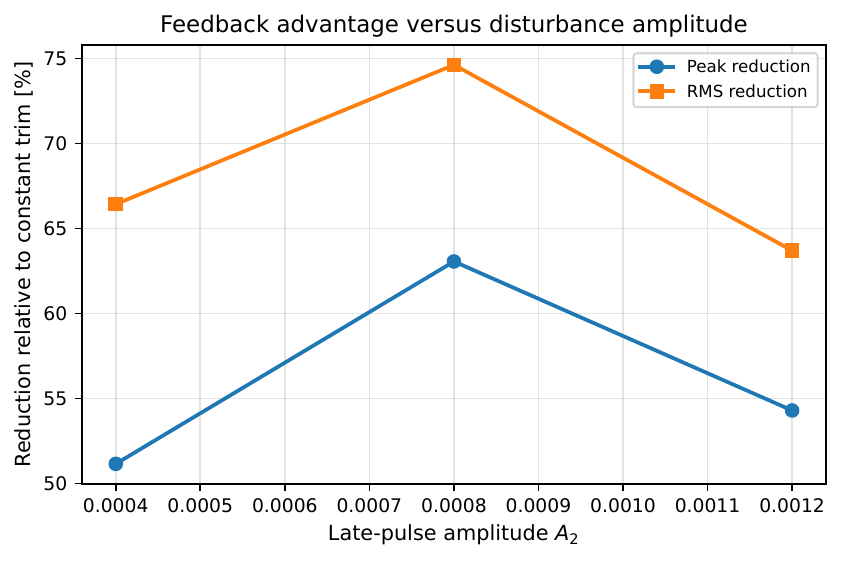}
\caption{SU2 Python-FSI disturbance-amplitude study. The bounded-PD baseline
retains a clear post-pulse advantage over constant trim across the three
perturbation amplitudes.}
\label{fig:fsi-robustness}
\end{figure}

Long-horizon 8000-step runs confirmed that the improvement was not confined to
the immediate pulse transient. For $A_2=8\times10^{-4}$, the RMS deviation
over 4--8~s fell from $7.49\times10^{-4}$ with constant trim to
$1.89\times10^{-4}$ with feedback, a 74.8\% reduction. The stronger
$A_2=1.2\times10^{-3}$ case produced essentially the same late-time reduction.
This result is reported as closed-loop FSI feasibility and disturbance
rejection, not as flutter suppression.

\begin{figure}[t]
\centering
\begin{subfigure}[t]{.49\linewidth}
\centering
\includegraphics[width=\linewidth,height=.34\textheight,keepaspectratio]{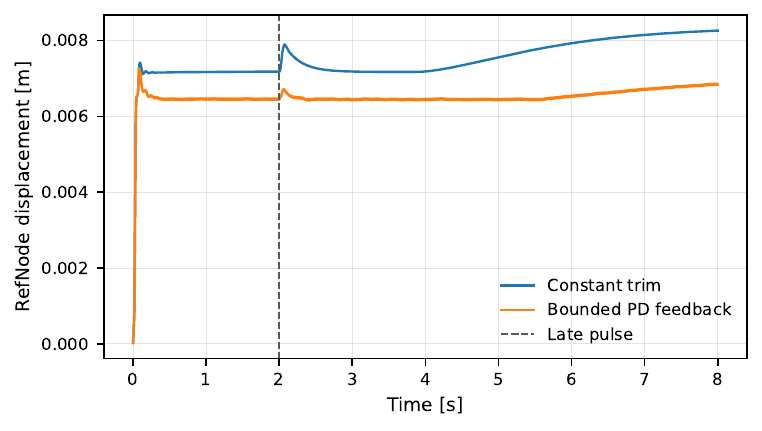}
\caption{Moderate late pulse, $A_2=8\times10^{-4}$.}
\end{subfigure}\hfill
\begin{subfigure}[t]{.49\linewidth}
\centering
\includegraphics[width=\linewidth,height=.34\textheight,keepaspectratio]{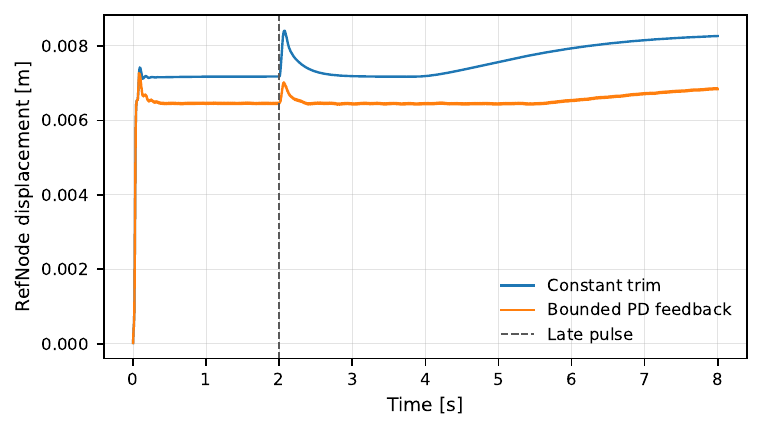}
\caption{Stronger late pulse, $A_2=1.2\times10^{-3}$.}
\end{subfigure}
\caption{Long-horizon SU2 Python-FSI validation. Constant trim maintains
numerical completion but exhibits a larger residual response; bounded
feedback provides sustained late-time attenuation.}
\label{fig:fsi-long}
\end{figure}

\FloatBarrier
\section{Discussion}
\subsection{What the primary benchmark demonstrates}
The primary benchmark establishes a complete data-to-control chain. HODMD
identifies the modal spectrum from measured structural channels; DMDc adds the
control input without using the true input matrix; the held-out chirp checks
that the resulting predictor is not merely a free-response fit; and MPC uses
that predictor to manage a hard flap constraint. The large-initial-condition
case is the most discriminating result because both controllers use the same
reduced model and the same input bound, yet only the constraint-aware MPC
avoids sustained clipping and keeps the true plant inside the stated linear
validity envelope.

\subsection{What the support studies add}
The support studies reduce three interpretation risks. SU2 shows that the
modal extraction is not confined to an analytical state-space generator.
OpenFOAM shows that an articulated boundary can deliver a nearly linear moment
response over the controller's deflection range while retaining important
unsteady frequency dependence. The Python-FSI study shows that a sensor,
controller and distributed structural load can be coupled around a two-way
multiphysics solve. None of these facts changes the level of the central
claim: the present MPC is validated against the separate typical-section true
plant, not directly against SU2 or OpenFOAM.

\subsection{Why MPC was not placed directly inside SU2}
A high-fidelity plant cannot normally be used as the internal prediction model
of an online MPC because every control update requires many forward
predictions. The coherent next step is therefore to excite the coupled plant,
identify a control-oriented ROM from SU2 data, and use that ROM inside the
optimizer while SU2 remains the external true plant. The current PD loop was a
necessary instrumentation and feasibility step, but it is retained only as a
baseline rather than promoted to the paper's principal controller.

\subsection{Limitations}
The typical-section benchmark is linear, two-dimensional and incompressible;
its post-critical records represent early-time exponential growth rather than
nonlinear limit-cycle oscillations. The flap is ideal and neglects actuator
rate limits, deadband and second-order dynamics. The OpenFOAM geometry is a
low-Reynolds-number demonstrator with an artificial hinge gap, not a
geometrically matched validation model. The SU2 modal cases do not establish a
flutter boundary, and the Python-FSI case uses structural-force actuation
rather than an aerodynamic flap. No experimental validation is included.
Finally, the five-state ROM is only mildly smaller than the six-state benchmark;
the computational advantage should become more relevant for CFD- or
finite-element-derived systems. These limitations motivate, rather than
invalidate, the planned transition to a ROM identified from a high-fidelity
coupled plant.

\FloatBarrier
\section{Conclusions}
A reproducible HODMD--DMDc--MPC workflow was developed for constrained
early-time flutter suppression. HODMD recovered the dominant stable and
unstable modes from four measured channels, and delay embedding preserved
near-critical growth-rate accuracy under 2\% RMS measurement noise. The modal
subspace was converted into a five-state real predictor by estimating the flap
input matrix from PRBS data; independent chirp and free-response tests yielded
normalized errors of $2.4\times10^{-3}$ and $7.7\times10^{-5}$, respectively.
At $1.05\,U_F$, the constrained MPC stabilized an $8^\circ$ initial pitch
perturbation with a $1^\circ$ flap bound, while saturated LQR reached the
$15^\circ$ linear-model cap and was terminated.

The CFD/FSI studies support, but do not enlarge, this central claim. SU2
confirmed dominant-frequency recovery from unsteady aerodynamic and
coupled-aeroelastic data. OpenFOAM showed approximately linear articulated
moment authority over $2^\circ$--$8^\circ$ and strong frequency dependence.
The SU2 Python-FSI baseline demonstrated sustained disturbance-rejection
benefits from bounded feedback over several pulse amplitudes and 8000-step
runs. The immediate next research step is a control-oriented ROM identified
from the coupled SU2 plant and used by constrained MPC in an external
high-fidelity closed loop.

\section*{Declaration of competing interest}
The author declares no known competing financial interests or personal
relationships that could have appeared to influence the work reported in this
paper.

\section*{Funding}
This research received no specific grant from funding agencies in the public,
commercial, or not-for-profit sectors.

\section*{CRediT authorship contribution statement}
\textbf{Carlos Domingo Mendez Gaona:} Conceptualization, Methodology,
Software, Validation, Formal analysis, Investigation, Visualization, Writing
-- original draft, Writing -- review and editing.

\section*{Data and code availability}
The numerically simulated benchmark datasets, identification scripts,
controller code, processed CFD/FSI data and figure-generation scripts are
provided as supplementary material. The raw high-fidelity solver outputs are
available from the author upon reasonable request because of their size.
\appendix
\section{Flap constants and rational-approximation coefficients}
\label{app:T}
For hinge location $c$ (in semichords), the flap constants of
Theodorsen's classical theory~\cite{theodorsen1935general} used in
Eq.~\eqref{eq:w34} are
\begin{equation}
T_{10} = \sqrt{1-c^2} + \arccos c , \qquad
T_{11} = \arccos c\,(1-2c) + \sqrt{1-c^2}\,(2-c),
\end{equation}
where $T_{11}$ multiplies the $\dot\beta$ downwash term activated only
in the actuator-dynamics extension. In addition to their classical
origin~\cite{theodorsen1935general,jones1940unsteady}, the numerical
adequacy of these constants and of the Wagner approximation coefficients ($a_1=0.165$,
$b_1=0.0455$, $a_2=0.335$, $b_2=0.300$) is established independently by
the verification checks of Section~2 and Fig.~\ref{fig:rfa}.


\end{document}